\documentclass{article}
\usepackage{array, amsmath, amssymb, graphicx, stmaryrd,a4wide}
\usepackage[percent]{overpic}
\usepackage{subfig}
 \usepackage{authblk}

\newcommand{\bea}[1][]{\begin{eqnarray}}
\newcommand{\tnorm}[1]{\interleave#1\interleave}
                                      
\newcommand{\eea}[1][]{\end{eqnarray}}

\newcommand{\bfn}{\boldsymbol n}
\newcommand{\bfx}{\boldsymbol x}

\newcommand{\bfp}{\boldsymbol p}
\newcommand{\bfI}{\boldsymbol I}
\newcommand{\bftau}{\boldsymbol \tau}
\newcommand{\bfP}{\boldsymbol P}
\newcommand{\bfv}{\boldsymbol v}
\newcommand{\bfu}{\boldsymbol u}
\newcommand{\bff}{\boldsymbol f}
\newcommand{\bfF}{\boldsymbol F}

\newcommand{\bfE}{\boldsymbol E}

\newcommand{\bfg}{\boldsymbol G}
\newcommand{\bfJ}{\boldsymbol J}

\newcommand{\bfkappa}{\boldsymbol\kappa}
\newcommand{\bfzero}{\boldsymbol 0}
\newcommand{\bfeps}{\boldsymbol\varepsilon}
\newcommand{\bfsig}{\boldsymbol\sigma}

\newcommand{\IR}{\mathbb{R}}
\newcommand{\ds}{~\mathrm{d}\Gamma}
\newcommand{\divergence}[1]{\nabla_\Gamma\cdot{#1}}
\newcommand{\divh}[1]{\nabla_{\Gamma_h}\cdot{#1}}

\begin{document}
\title{Augmented Lagrangian and Galerkin least squares methods for membrane contact} 
\author[$\dagger$]{Erik Burman}
\author[$\ddag$]{Peter Hansbo}
\author[$\star$]{Mats G.\ Larson}
\affil[$\dagger$]{\small Department of Mathematics, University College London, London, UK--WC1E  6BT, United Kingdom}
\affil[$\ddag$]{\small Department of Mechanical Engineering, J\"onk\"oping University, SE-55111 J\"onk\"oping, Sweden}
\affil[$\star$]{\small Department of Mathematics and Mathematical Statistics, Ume{\aa} University, SE-901 87 Ume{\aa}, Sweden}
\maketitle

\abstract{In this paper, we propose a stabilised finite element
method for the numerical
solution of contact between a small deformation elastic membrane and a rigid obstacle. We limit ourselves to 
friction--free contact, but the formulation is readily extendable to 
more complex situations.}

\section{Introduction}
Finite element solvers for contact problems typically employ either Lagrange multipliers or the penalty method for the mortaring on contact zones. The penalty method
is simple to implement and robust but inconsistent and thus requires a large penalty parameter to ensure non-penetration, which leads to ill conditioning and 
possibly instability (depending on how the penalty is imposed). Lagrange multipliers, on the other hand, require careful matching of the spaces for the 
primal variable and multiplier. Stabilised multiplier methods have been proposed as a remedy by, e.g., Heintz and Hansbo \cite{HeHa06}, Hild and Renard \cite{HiRe10}, Oliver et al. \cite{OlHa09}. These methods are often of Galerkin/Least Squares (GLS) type, where a penalty is placed on the deviation between the multiplier and the contact force (derived 
from the primal variable), an approach first proposed by Barbosa and Hughes \cite{BaHu91} for the linear boundary multiplier method.

A combination of the multiplier and penalty approaches yields the augmented Lagrangian method, cf. Alart and Curnier \cite{AlCu91}.
Provided the penalty is not too strong this is expected 
to improve the conditioning as well as improving the control of the constraint compared to the case where only the 
multiplier is used. For early work on augmented Lagrangian methods in computational methods 
for partial differential equations we refer to Fortin and Glowinski \cite{FG83}, The augmented Lagrangian approach was recently used by Chouly and Hild \cite{ChHi12}
to eliminate the multiplier and arrive at a Nitsche type method for contact, and their work was adapted to the case of the
the obstacle problem by Burman, Hansbo, and Larson \cite{BuHaLaSt17}. In this work we further develop the idea to handle the case of friction free contact between curved membranes and rigid obstacles, using tangential differential calculus for the membrane model \cite{HaLa14, HaLaLa15}.

The rest of the paper is organised as
follows. In Section \ref{membrane} we recall the membrane model from Hansbo and Larson \cite{HaLa14}, in Section \ref{augment} we describe the continuous and discrete versions of the proposed augmented Lagrangian method, and in Section \ref{GLS} we derive our GLS method. 
In Section \ref{numex}, we present some numerical results, and, finally, in Section \ref{conclusions} we give some concluding remarks.

\section{The membrane problem}
\label{membrane}
\subsection{Basic notation}

Let $\Gamma$ be a smooth two-dimensional surface embedded in $\IR^3$, with outward pointing normal $\bfn$. We shall here for simplicity assume that the surface is closed, but this is not a requirement for the following, boundary conditions can be applied as discussed in Hansbo and Larson \cite{HaLa14}.
If we denote the signed distance function relative to $\Gamma$ by $d(\bfx)$, for $\bfx\in \Bbb{R}^3$, fulfilling $\nabla d(\bfx) = \bfn(\bfx)$ if $\bfx\in \Gamma$, we can define the domain occupied by the membrane by
\[
\Omega_t = \{\bfx\in \Bbb{R}^3: \vert d(\bfx) \vert < t/2\},
\]
where $t$ is the thickness of the membrane. The closest point projection $\bfp:\Omega_t \rightarrow \Gamma$
is given by
\[
\bfp(\bfx) = \bfx -d(\bfx)\nabla d(\bfx) ,
\]
the Jacobian matrix of which is
\[
\nabla \bfp = \bfI -d ( \nabla\otimes\nabla d ) -\nabla d\otimes\nabla d
\]
where $\bfI$ is the identity and $\otimes$ denotes the exterior product 
$(a\otimes b)_{ij} = a_i b_j$ for vectors $a$ and $b$ in $\IR^3$.
The corresponding linear projector
$\bfP_\Gamma = \bfP_\Gamma(\bfx)$, onto the tangent plane of $\Gamma$
at $\bfx\in\Gamma$, is given by
\[
\bfP_\Gamma := \bfI -\bfn\otimes\bfn ,
\]
and we can use it to define the surface gradient
$\nabla_\Gamma$ as
\begin{equation}
\nabla_\Gamma := \bfP_\Gamma \nabla .
\end{equation}
The surface gradient thus has three
components, which we shall denote by
\[
\nabla_\Gamma =: \left[\begin{array}{>{\displaystyle}c}
\frac{\partial}{\partial x^\Gamma} \\[3mm]
\frac{\partial}{\partial y^\Gamma} \\[3mm]
\frac{\partial}{\partial z^\Gamma} \
\end{array}\right] .
\]
For a vector valued function $\bfv(\bfx)$, we define the
tangential Jacobian matrix as the transpose of the outer product of $\nabla_\Gamma$ and $\bfv$,
\[
\left(\nabla_\Gamma\otimes\bfv\right)^{\text{T}} :=\left[\begin{array}{>{\displaystyle}c>{\displaystyle}c>{\displaystyle}c}
\frac{\partial v_1}{\partial x^\Gamma} &\frac{\partial v_1}{\partial y^\Gamma} & \frac{\partial v_1}{\partial z^\Gamma} \\[3mm]
\frac{\partial v_2}{\partial x^\Gamma} &\frac{\partial v_2}{\partial y^\Gamma} & \frac{\partial v_2}{\partial z^\Gamma} \\[3mm]
\frac{\partial v_3}{\partial x^\Gamma} &\frac{\partial v_3}{\partial y^\Gamma} & \frac{\partial v_3}{\partial z^\Gamma}
\end{array}\right]
\]
and the surface divergence $\nabla_{\Gamma}\cdot\bfv := \text{tr}(\nabla_\Gamma\otimes\bfv)$.
%
%

\subsection{The surface strain and stress tensors}

To obtain an in-plane strain tensor we need to use the projection twice to define
\[
\bfeps_{\Gamma}(\bfu) := \bfP_\Gamma\bfeps(\bfu)\bfP_\Gamma ,
\]
which lacks all out-of-plane strain components. In other words both the rows and the columns of 
$\bfeps_{\Gamma}(\bfu)$ are tangent vectors so that $\bfeps_{\Gamma}(\bfu) \cdot \bfn = 
\bfn \cdot \bfeps_{\Gamma}(\bfu) = \bfzero$.
For a membrane, where plane stress is assumed,
this strain tensor can still be used, since out-of-plane strains do not contribute to the strain energy.
However, the tensor $\bfeps_{\Gamma}$ is rather cumbersome to use directly in a numerical implementation;
it is easier to work with the symmetric part of the surface Jacobian 
\[
\bfE_{\Gamma}(\bfu) := \frac12\left(\nabla_\Gamma\otimes \bfu + (\nabla_\Gamma\otimes\bfu)^{\rm T}\right),
\]
which can
be established directly using tangential derivatives. For this reason, we use the fact that
$\bfn\cdot\bfE_{\Gamma}(\bfu)\cdot\bfn=0$ to obtain the following relation:
\[
\bfeps_{\Gamma}(\bfu) 
= \bfP_\Gamma \bfE_\Gamma \bfP_\Gamma 
= \bfE_{\Gamma}(\bfu) - \left((\bfE_{\Gamma}(\bfu)\cdot\bfn)\otimes\bfn + \bfn\otimes (\bfE_{\Gamma}(\bfu)\cdot\bfn)\right) ,
\]
so that, using dyadic double-dot product,
\[
\bfsig:\bfu\otimes\bfv = (\bfsig\cdot\bfu)\cdot\bfv,
\]
where $\bfsig$ is a tensor and $\bfu$, $\bfv$ are vectors, we arrive at
\begin{equation}\label{doubledot}
\bfeps_{\Gamma}(\bfu):\bfeps_{\Gamma}(\bfv) = \bfE_{\Gamma}(\bfu):\bfE_{\Gamma}(\bfv) - 2(\bfE_{\Gamma}(\bfu)\cdot\bfn)\cdot(\bfE_{\Gamma}(\bfv)\cdot\bfn),
\end{equation}
which will be used in the finite element implementation below. 

We shall assume an isotropic stress--strain relation,
\[
\bfsig(\bfu) = 2\mu \bfeps(\bfu) + \lambda \nabla\cdot\bfu\, \bfI ,
\]
where $\bfsig$ is the stress tensor and $\bfI$ is the identity tensor.
The Lam\'e parameters $\lambda$ and $\mu$ are related to Young's modulus $E$
and Poisson's ratio $\nu$ via
\[
\mu = \frac{E}{2(1+\nu)},\quad \lambda =\frac{ E\nu}{(1+\nu)(1-2\nu)} .
\]
For
the in-plane  stress tensor we assume
\begin{equation}\label{constitutive}
\bfsig_\Gamma(\bfu):=  2\mu \bfeps_\Gamma(\bfu) + {\lambda_0}\nabla_\Gamma\cdot\bfu\, \bfP_\Gamma ,
\end{equation}
where
\[
{\lambda_0}:= \frac{2\lambda\mu}{\lambda+2\mu} = \frac{E\nu}{1-\nu^2}.
\]
is the Lam\'e parameter in plane stress conditions. This assumption is consistent with the membrane model of Ciarlet and Sanchez-Palencia \cite{CiSa93}, as shown by Delfour and Zol\'esio \cite{DeZo97}.
We remark that out-of-plane components of the contraction between stress and strain will not contribute anything to the strain energy functional underlying the finite element method, since both tensors are in-plane. Thus, the only difference 
between plane stress and plane strain in a curved membrane (as concerns strain energy) lies in the distinction between $\lambda$ and $\lambda_0$, as in the two-dimensional case.

\section{Augmented Lagrangian formulation of the membrane contact problem}\label{augment}

The equilibrium equation for the membrane can be written
 \begin{equation}\label{equil}
 -\nabla_\Gamma\cdot \bfsig_{\Gamma}(\bfu) = \bff\quad\text{in}\; \Gamma ,
 \end{equation}
where the matrix divergence is defined by taking the vector surface divergence of each row of $\bfsig_{\Gamma}$, cf. Hansbo and Larson \cite{HaLa14}.
Note that here $\bff\in L_2(\Gamma)$ is proportional to $t^{-1}$ (so that $\bff$ has units force per unit volume). Equation (\ref{equil}), together with the constitutive law (\ref{constitutive}) defines the 
 differential equations of linear elasticity in general on surfaces.

Our model problem of friction free 
contact between the membrane and a rigid obstacle thus takes the form
\begin{align}
\bfP_\Gamma (\bff + \divergence{\bfsig_{\Gamma}})  = {\boldsymbol 0} &{} \quad \text{in}\; \Gamma ,\\ \label{strongform1}
\bfn\cdot(\bff + \divergence{\bfsig_{\Gamma}})  \geq  0 &{}  \quad \text{in}\; \Gamma ,\\ 
\bfsig_{\Gamma} = 2\mu\bfeps_\Gamma(\bfu)+\lambda_0\nabla_\Gamma\cdot\bfu\,\bfP_\Gamma  &{} \quad \text{in}\; \Gamma , \\ (u_{n}-g) \leq 0 ,\;
(u_{n}-g)\,\bfn\cdot(\bff + \divergence{\bfsig_{\Gamma}}) =  0 &{}\quad \text{in}\; \Gamma ,
\end{align}
with 
$u_n:=\bfu\cdot\bfn$ and $g$ denotes the normal distance from the membrane to the obstacle before deformation.

In order to define the augmented Lagrangian method, we first introduce a Lagrange multiplier $p$ such that
\begin{align}
-\divergence{\bfsig_{\Gamma}} - p\,\bfn = \bff   &{} \quad \text{in}\; \Gamma ,\\ \label{kuhn1}
u_n-g  \leq  0 &{}  \quad \text{in}\; \Gamma ,\\ \label{kuhn2}
p \leq 0  &{} \quad \text{in}\; \Gamma , \\ \label{kuhn3}
(u_{n}-g)\,p =  0 &{}\quad \text{in}\; \Gamma , 
\end{align}
Note that the multiplier here has the interpretation as the out-of-balance normal force per unit volume $p := -\bfn\cdot(\bff + \divergence{\bfsig_{\Gamma}})$ stopping the membrane from penetrating the rigid object.

We can now replace the Kuhn--Tucker conditions (\ref{kuhn1})--(\ref{kuhn3}), using the notation
\begin{equation}\label{positivepart}
[a]_+ := \left\{\begin{array}{l} a\quad\text{if $a>0$},\\ 0\quad \text{if $a\leq 0$,}\end{array}\right. ,
\end{equation}
 by the equivalent statement
\begin{equation}\label{lambdadef}
p = -\frac{1}{\gamma}[u_{n}-g -\gamma p]_+
\end{equation}
with $\gamma$ a positive number, cf . Chouly and Hild \cite[Prop. 2.1]{ChHi12}. Note that for dimensional reasons $\gamma$ must be proportional to the thickness squared and inversely  proportional to the Lam\'e parameters. 

Defining the natural function space for the displacements as
\[
V = \{\bfv  : v_n \in L_2(\Gamma )\quad\text{and}\quad \bfv-\bfn\, v_n =:\bfv_t \in [H^1(\Gamma)]^2\},
\]
cf. Ciarlet et al. \cite{CiSa93,CiLo96}, and for the multipliers as
\begin{equation}\label{funcV}
 Q = L_2(\Gamma),
\end{equation}
and seeking $(\bfu,p)\in V\times Q$ we have by Green's theorem on surfaces (cf. Gurtin and Murdoch\cite{GuMo75}), with
\[
 L(\bfv) := \int_{\Gamma} \bff\cdot\bfv\, {\rm d}\Gamma ,
\]
and
\begin{align*}
a_{\Gamma}(\bfu,\bfv)  := &{} \int_{\Gamma} \bfsig_{\Gamma}(\bfu):\bfeps_{\Gamma}(\bfv)\, {\rm d}\Gamma \\
=&{}\int_{\Gamma}2\mu\bfeps_\Gamma(\bfu): \bfeps_\Gamma(\bfv){\,\rm d}\Gamma+ \int_{\Gamma}{\lambda_0}\, \nabla_\Gamma\cdot\bfu,\nabla_\Gamma\cdot\bfv{\,\rm d}\Gamma\\
= {}& \int_{\Gamma}2\mu\bfE_\Gamma(\bfu): \bfE_\Gamma(\bfv){\,\rm d}\Gamma- \int_{\Gamma}4\mu\bfE_\Gamma(\bfu)\cdot\bfn, \bfE_\Gamma(\bfv)\cdot\bfn{\,\rm d}\Gamma \\
 & +  \int_{\Gamma}{\lambda_0}\, \nabla_\Gamma\cdot\bfu,\nabla_\Gamma\cdot\bfv{\,\rm d}\Gamma,
\end{align*}
that
\[
a_{\Gamma}(\bfu,\bfv) -\int_{\Gamma} p \, v_n\, \ds = L(\bfv)
\]
where $\bfv\in V$. Following Chouly and Hild \cite{ChHi12} we write $v_n = v_n + \gamma q-\gamma q$
for an arbitrary function $ q\in Q$, so that we may write
\[
a_{\Gamma}(\bfu,\bfv) -\int_{\Gamma} p \, (v_n-\gamma q) \, \ds -\int_{\Gamma}\gamma p \,  q \, \ds = L(\bfv) .
\]
Replacing $p $ in the first integral by the expression in (\ref{lambdadef}) we finally obtain the problem of finding $(\bfu,p )\in V\times Q$ such that
\begin{equation}\label{aug}
a_{\Gamma}(\bfu,\bfv) +\int_{\Gamma} \frac{1}{\gamma} [u_{n}-g -\gamma p ]_+ (v_n-\gamma q) \, \ds -\int_{\Gamma}\gamma p \,  q \, \ds = L(\bfv) 
\end{equation}
for all $(\bfv, q)\in V\times Q$.
This problem is related to seeking stationary points to the functional
\begin{equation}\label{augmin}
\Pi(\bfu,p ) := a_{\Gamma}(\bfu,\bfu) - L(\bfu) + \int_{\Gamma}\frac{1}{2\gamma}\left[u_n-g-\gamma p \right]_+^2\ds -\int_{\Gamma}\frac{\gamma}{2} p ^2\ds,
\end{equation}
see, e.g., Alart and Curnier \cite{AlCu91}.  
The formulation (\ref{augmin}) constitutes the starting point for our
finite element approximation.

\subsection{The finite element method}

Let $\mathcal{T}_h:=\{T\}$ be a conforming, shape regular triangulation of {$\Gamma$} using a parametric map of a certain polynomial degree from reference triangles, resulting in a discrete surface $\Gamma_h$ constructed as follows.
We wish to define a map $\bfF : (\xi,\eta)\rightarrow (x,y,z)$ from a reference triangle $\hat T$ defined in a local coordinate system $(\xi, \eta)$ to $T$, for all $T$. To this end, we write
$\bfx_\Gamma = \bfx_\Gamma(\xi, \eta)$, where $\bfx_\Gamma=(x_\Gamma,y_\Gamma,z_\Gamma)$
are the physical coordinates on $\Gamma$.
For any given parametrization, we can extend it outside the surface by
defining
\[
\bfx(\xi,\eta,\zeta) = \bfx_\Gamma(\xi, \eta)+\zeta\,\bfn(\xi,\eta)
\]
where $\bfn$ is the normal and $-t/2 \leq \zeta \leq t/2$. 
We next consider an elementwise parametrization of the surface as
\[
\bfx_{\Gamma_h}(\xi,\eta)= \sum_i\bfx_\Gamma^i\psi_i(\xi,\eta)
\]
where $\bfx_\Gamma^i$ are the coordinates of the nodes, assumed located on $\Gamma$, and $\psi_i(\xi,\eta)$ are finite element shape functions of a certain degree on the reference element, and extend this approximation outside the surface so that
\begin{equation}\label{approx}
\bfx (\xi,\eta,\zeta) \approx \bfx^h(\xi,\eta,\zeta):= \bfx_{\Gamma_h}(\xi,\eta)+\zeta\,\bfn^h(\xi,\eta)
\end{equation}
 where
\begin{equation}\label{normal}
\bfn^h = \frac{\displaystyle\frac{\partial \bfx_{\Gamma_h}}{\partial \xi}\times\frac{\partial \bfx_{\Gamma_h}}{\partial \eta}}{\displaystyle\left\vert \frac{\partial \bfx_{\Gamma_h}}{\partial \xi}\times\frac{\partial \bfx_{\Gamma_h}}{\partial \eta}\right\vert}.
\end{equation}
This gives us the exact normal vector to the discrete surface, which ensures that the correct rigid body motions are reproduced in the discrete model, i.e., that
\begin{equation}
\text{ker}\,\bfeps_{\Gamma_h}(\bfv) := \{\bfv\in [H^1(\Gamma_h)]^3:\; \bfeps_{\Gamma_h}(\bfv) = {\boldsymbol 0}\;\text{and}\; \bfv\cdot\bfn^h=0\}
\end{equation}
is finite dimensional and consists only of rigid body rotations $\bfv$.

For the approximation of the displacement, we use a constant extension,
\begin{equation}
\bfu \approx \bfu^h = \sum_i\bfu_i\varphi_i(\xi,\eta)
\end{equation}
where $\bfu_i$ are the nodal displacements, and $\varphi_i$ are shape functions, not necessarily of the same degree as the $\psi_i$. Note that only the in-plane variation of the approximate solution will matter since we are looking at in-plane stresses and strains.
We employ the usual finite element approximation of the physical derivatives of the chosen basis $\{\varphi_i\}$ on the surface, at $(\xi,\eta,0)$, as
\begin{equation}\label{jacobiinv}
\left[\begin{array}{>{\displaystyle}c}
\frac{\partial \varphi_i}{\partial x}\\[2mm]
\frac{\partial \varphi_i}{\partial y}\\[2mm]
\frac{\partial \varphi_i}{\partial z}\end{array}\right] = \bfJ^{-1}(\xi,\eta,0) \left[\begin{array}{>{\displaystyle}c}
\frac{\partial \varphi_i}{\partial \xi}\\[2mm]
\frac{\partial \varphi_i}{\partial \eta}\\[2mm]
0\end{array}\right] ,
\end{equation}
where
\[
\bfJ(\xi,\eta,0) := \left[\begin{array}{>{\displaystyle}c>{\displaystyle}c>{\displaystyle}c}
\frac{\partial x_{\Gamma_h}}{\partial \xi} & \frac{\partial y_{\Gamma_h}}{\partial \xi} & \frac{\partial z_{\Gamma_h}}{\partial \xi}\\[3mm]
\frac{\partial x_{\Gamma_h}}{\partial \eta} & \frac{\partial y_{\Gamma_h}}{\partial \eta} & \frac{\partial z_{\Gamma_h}}{\partial \eta}\\[3mm]
n^h_y & n^h_y & n^h_z\end{array}\right] ,
\]
cf. Hansbo and Larson \cite{HaLa14}.

\section{Galerkin least squares method}\label{GLS}
\subsection{Formulation}

As a first attempt at a finite element method for the membrane contact problem we can introduce finite element spaces constructed from the basis previously discussed by defining
\begin{equation}\label{spacevA}
W^h_k:= \{ v: {v\vert_T \circ\bfF\in P^k(\hat T),\; \forall T\in\mathcal{T}_h};\; v \in C^0(\Gamma_h)\},
\end{equation}
and
\[
Q^h =\{ v: {v\vert_T \circ\bfF\in P^l(\hat T),\; \forall T\in\mathcal{T}_h}\},
\]
$l\leq k$,
and a tentative finite element method to be reads: Find $(\bfu^h,p^h)\in V^h\times Q^h$, where $V^h := [W^h_1]^3$, such that
\begin{equation}\label{augdisc}
a_{\Gamma_h}(\bfu^h,\bfv) +\int_{\Gamma_h} \frac{1}{\gamma} [u^h_{n}-g -\gamma p^h ]_+ (v_n-\gamma q)  \,\text{d}\Gamma_h -\int_{\Gamma_h}\gamma p^h \,  q \,\text{d}\Gamma_h = L(\bfv) 
\end{equation} 
for all $(\bfv, q)\in V^h\times Q^h$
where
\[
a_{\Gamma_h}(\bfu,\bfv)   := \int_{\Gamma_h}\bfsig_{\Gamma_h}(\bfu): \bfeps_{\Gamma_h}(\bfv))\,\text{d}\Gamma_h, \qquad L(\bfv) :=\int_{\Gamma_h}\bff^e\cdot\bfv\,\text{d}\Gamma_h
\]
where $\bff^e$ denotes an extension of $\bff$ from $\Gamma$ to $\Gamma_h$. Here and below we write $u_n^h := \bfu^h\cdot\bfn^h$ for $\bfu^h\in V^h$. The discrete geometry $\Gamma_h$ is constructed from a finite element interpolation of 
the exact geometry in $[W^h_m]^3$, with $m$ not necessarily equal to $k$. 

Clearly, not all combinations of discrete spaces are stable;
we therefore apply a GLS stabilization method which also allows for the elimination of the pressure variable.
To this end, we formally replace $p$ in (\ref{augmin}) by
$-\bfn\cdot(\bff + \divergence{\bfsig_{\Gamma}(\bfu)})$ and
seek $\bfu^h\in V_h$ such that
\begin{equation}
\bfu^h = \arg\min_{\bfv\in V_h} \mathfrak{F}_h(\bfv)
\end{equation}
where
\begin{align}\nonumber
\mathfrak{F}_h(\bfv) := {}& \frac12 a_{\Gamma_h}(\bfv,\bfv)  
+ \int_{\Gamma_h}\frac{1}{2\gamma}\left[ v_n-g +\gamma\bfn^h\cdot(\bff + \divh{\bfsig_{\Gamma_h}(\bfv)})\right]_+^2\, \text{d}\Gamma_h  \\
{}& -\frac12\int_{\Gamma_h} \gamma (\bfn^h\cdot(\bff + \divh{\bfsig_{\Gamma_h}(\bfv)}))^2 \, \text{d}\Gamma_h
-L(\bfv)  .\label{perturbed}
\end{align}
The
Euler--Lagrange equations corresponding to (\ref{perturbed}) take the form: Find $\bfu^h \in V_h$ such that
\begin{equation}\label{FEM}
a_{\Gamma_h}(\bfu^h,\bfv) + b_{\Gamma_h}(\bfu^h,g,\bff;\bfv) =L(\bfv) \quad \forall \bfv
\in V_h 
\end{equation}
where 
\begin{align}
b_{\Gamma_h}(\bfu,g,\bff;\bfv)
:= &{}
\int_{\Gamma_h}
 \gamma^{-1}[u_n- g+\gamma \bfn^h\cdot( \bff+ \divh{\bfsig_{\Gamma_h}(\bfu)})]_+ (v_n + \gamma \bfn^h \cdot (\divh{\bfsig_{\Gamma_h}(\bfv)})) \,\text{d}{\Gamma_h} 
 \\ 
 & 
 -   \int_{\Gamma_h}\gamma \bfn^h\cdot( \bff+ \divh{\bfsig_{\Gamma_h}(\bfu)}) 
 \,  \bfn^h\cdot (\divh{\bfsig_{\Gamma_h}(\bfv))}\,\text{d}{\Gamma_h} .\label{stab_form0}
\end{align}

Next we have the identity
\begin{equation}\label{divsigma}
\bfn\cdot(\divergence{\bfsig_{\Gamma}(\bfu))} = -\bfsig_{\Gamma}(\bfu):\bfkappa 
\end{equation}
where $\bfkappa := \nabla \otimes \bfn$ is the curvature tensor (negative Weingarten map) and $\bfsig_{\Gamma}(\bfu):\bfkappa$ denotes the Frobenius inner product. To verify (\ref{divsigma}) we note that multiplying with a test function $\varphi$ and using integration by parts we obtain, for tangential $\bfsig$,
\begin{align}
\int_\Gamma \bfn \cdot (\divergence \bfsig) \varphi \,\text{d}{\Gamma}
&=
-\int_\Gamma \bfsig : \nabla (\varphi \bfn)\,\text{d}{\Gamma}
=
- \int_\Gamma \bfsig :(\nabla \varphi) \otimes \bfn\,\text{d}{\Gamma}
- \int_\Gamma \bfsig :\bfkappa \varphi\,\text{d}{\Gamma}
=
- \int_\Gamma \bfsig :\bfkappa \varphi\,\text{d}{\Gamma}
\end{align}
where we used the fact that $\bfsig$ is tangential to conclude that 
$\bfsig :(\nabla \varphi) \otimes \bfn = (\nabla \varphi) \cdot \bfsig \cdot \bfn =0$. Therefore we conclude that 
\begin{equation}
\int_\Gamma ( \bfn \cdot (\divergence \bfsig) + \bfsig :\bfkappa) \varphi\,\text{d}{\Gamma} = 0  
\end{equation}
and thus (\ref{divsigma}) holds.
Furthermore, we note that $\bfP:\bfkappa =\text{tr }\bfkappa$ and $\bfeps_{\Gamma}:\bfkappa = \bfE_\Gamma:\bfkappa$ so that
\[
\bfsig_{\Gamma}(\bfu):\bfkappa = 2\mu  \bfE_\Gamma(\bfu):\bfkappa + \lambda_0 \nabla_\Gamma\cdot\bfu\;\text{tr }\bfkappa
\]
which together with (\ref{divsigma}) can be used to simplify the implementation.


In the discrete version of (\ref{divsigma}) we need the approximate curvature tensor $\bfkappa^h :=\nabla \otimes\bfn^h$. Since we have no explicit expression for
$\bfn^h$ in physical coordinates, 
we must use some tools from classical differential geometry to compute $\bfkappa^h$.
To this end, we define the matrix representations of the first and second fundamental forms as follows:
\[
\bfg_1:=\left[\begin{array}{>{\displaystyle}c>{\displaystyle}c}
\frac{\partial \bfx_\Sigma^h}{\partial \xi}\cdot\frac{\partial \bfx_\Sigma^h}{\partial \xi} & \frac{\partial \bfx_\Sigma^h}{\partial \xi}\cdot\frac{\partial \bfx_\Sigma^h}{\partial\eta}\\[3mm]
\frac{\partial \bfx_\Sigma^h}{\partial \xi}\cdot\frac{\partial \bfx_\Sigma^h}{\partial\eta} & \frac{\partial \bfx_\Sigma^h}{\partial \eta}\cdot\frac{\partial \bfx_\Sigma^h}{\partial \eta}
\end{array}\right] 
\]
\[
\bfg_2:=\left[\begin{array}{>{\displaystyle}c>{\displaystyle}c}
\bfn^h\cdot\frac{\partial^2 \bfx_\Sigma^h}{\partial \xi^2} & \bfn^h\cdot\frac{\partial^2 \bfx_\Sigma^h}{\partial \xi\partial\eta}\\[3mm]
\bfn^h\cdot\frac{\partial^2 \bfx_\Sigma^h}{\partial \xi\partial\eta} & \bfn^h\cdot\frac{\partial^2 \bfx_\Sigma^h}{\partial \eta^2}
\end{array}\right] 
\]
Through the Weingarten equations, cf. Kreyszig \cite{Kr91}, we may then compute the $(\xi,\eta)$--derivatives of the components of $\bfn^h$:
\[
\left[\begin{array}{>{\displaystyle}c>{\displaystyle}c>{\displaystyle}c} \frac{\partial n_x^h}{\partial \xi} &  \frac{\partial n_y^h}{\partial \xi}& \frac{\partial n_z^h}{\partial \xi}\\[3mm]
\frac{\partial n_x^h}{\partial \eta} & \frac{\partial n_y^h}{\partial \eta} & \frac{\partial n_z^h}{\partial \eta}\end{array}\right] =-\bfg_1^{-1}\bfg_2 \left[\begin{array}{>{\displaystyle}c>{\displaystyle}c>{\displaystyle}c} \frac{\partial x_\Sigma^h}{\partial \xi} & \frac{\partial y_\Sigma^h}{\partial \xi} & \frac{\partial z_\Sigma^h}{\partial \xi} \\[3mm]
\frac{\partial x_\Sigma^h}{\partial \eta} & \frac{\partial y_\Sigma^h}{\partial \eta} & \frac{\partial z_\Sigma^h}{\partial \eta}\end{array}\right].
\]
The physical derivatives of $\bfn^h$, defining the curvature tensor $\bfkappa^h$ in 3D, are then found by use of the Jacobian matrix $\bfJ$ analogously to (\ref{jacobiinv}).

Denoting
\[
\sigma_\kappa(\bfu) := -2\mu  \bfE_{\Gamma_h}(\bfu):\bfkappa^h - \lambda_0 \nabla_{\Gamma_h}\cdot\bfu\;\text{tr }\bfkappa^h ,\qquad f^h_n:=\bfn^h\cdot\bff,
\]
our bilinear form $b_{\Gamma_h}(\cdot,\cdot,\cdot;\cdot)$ takes the form
\begin{equation}
b_{\Gamma_h}(\bfu^h,g,\bff;\bfv^h):= \int_{\Gamma_h}
 \left(\gamma^{-1}[u^h_n- g+\gamma (f_n^h+ \sigma_\kappa(\bfu^h))]_+( v^h_n + \gamma \sigma_\kappa(\bfv^h))\,
 -   \gamma (f^h_n+\sigma_\kappa(\bfu^h))  \sigma_\kappa(\bfv^h)\right)\,\text{d}{\Gamma_h}. \label{stab_form1}
\end{equation}

\subsection{On the stability of the method}

We introduce the discrete linear operator
$P_\gamma(\bfv^h) := v^h_n + \gamma \sigma_\kappa(\bfv^h)$, and use the notation $\psi := g-\gamma f^h_n$ so that
\begin{equation}
b_{\Gamma_h}(\bfu^h,g,\bff;\bfv^h)= \int_{\Gamma_h}
 \left(\gamma^{-1}[P_\gamma(\bfu^h)- \psi]_+P_\gamma(\bfv^h)\,
 -   \gamma (f^h_n+\sigma_\kappa(\bfu^h))  \sigma_\kappa(\bfv^h)\right)\,\text{d}{\Gamma_h}, \label{stab_form3}
\end{equation}
and we note that $b_{\Gamma_h}$ can be interpreted as a nonlinear
penalty term, consistent on $\Gamma_h$, for the imposition of the
contact condition. Taking $\bfv_h=\bfu_h$ in \eqref{stab_form3} leads to 
\begin{align}\nonumber
b_{\Gamma_h}(\bfu^h,g,\bff;\bfu^h)= {}& \int_{\Gamma_h}
 \left(\gamma^{-1}[P_\gamma(\bfu^h)- \psi]_+P_\gamma(\bfu^h)\,
 -   \gamma (f^h_n+\sigma_\kappa(\bfu^h))
 \sigma_\kappa(\bfu^h)\right)\,\text{d}{\Gamma_h} \\
= {}&\gamma^{-1}\|[P_\gamma(\bfu^h)- \psi]_+\|^2_{L_2(\Gamma_h)} -  \gamma\| \sigma_\kappa(\bfu^h)\|_{L_2(\Gamma_h)}^2+ \int_{\Gamma_h}
 \left(\gamma^{-1}[P_\gamma(\bfu^h)- \psi]_+\psi \, -   \gamma  f^h_n
   \sigma_\kappa(\bfu^h)\right) \,\text{d}{\Gamma_h}.\label{bstab}
\end{align}
Using this relation in the formulation \eqref{FEM} leads to the
equality
\begin{equation}\begin{array}{>{\displaystyle}c}\label{stab_rel_1}
a_{\Gamma_h}(\bfu^h,\bfu^h)+\gamma^{-1}\|[P_\gamma(\bfu^h)- \psi]_+\|^2_{L_2(\Gamma_h)}- \gamma\| \sigma_\kappa(\bfu^h)\|_{L_2(\Gamma_h)}^2 
= \\
L(\bfu^h)+\int_{\Gamma_h}\left(\gamma  f^h_n
   \sigma_\kappa(\bfu^h)-
 \gamma^{-1}[P_\gamma(\bfu^h)- \psi]_+\psi    \right) \,\text{d}{\Gamma_h}. \end{array}\end{equation}
For $\gamma$ small enough, we have that
\begin{equation}\label{stabdelta}
a_{\Gamma_h}(\bfu^h,\bfu^h)- \gamma\| \sigma_\kappa(\bfu^h)\|_{L_2(\Gamma_h)}^2 \geq C_\gamma a_{\Gamma_h}(\bfu^h,\bfu^h)
\end{equation}
with $C_\gamma$ a constant independent of the meshsize. To show this, note that 
\[
\vert\sigma_\kappa(\bfu)\vert \leq \vert \bfsig_{\Gamma_h}(\bfu):\bfkappa^h\vert \leq  \vert \bfsig_{\Gamma_h}(\bfu)\vert\, \vert\bfkappa^h\vert
\]
where $\vert\bftau\vert$ denotes the Frobenius norm of a matrix
$\bftau$ and $\vert\bftau\vert_\text{max}$ its maximum over $\Gamma_h$. Thus
\begin{align}\nonumber
a_{\Gamma_h}(\bfu^h,\bfu^h)- \gamma\| \sigma_\kappa(\bfu^h)\|_{L_2(\Gamma_h)}^2 & \geq a_{\Gamma_h}(\bfu^h,\bfu^h)- \gamma\vert\bfkappa^h\vert_\text{max}^2\|  \bfsig_{\Gamma_h}(\bfu^h)\|_{L_2(\Gamma_h)}^2 \\  &\geq  a_{\Gamma_h}(\bfu^h,\bfu^h)\left(1-\gamma\vert\bfkappa^h\vert_\text{max}^2(2\mu+2\lambda_0)\right),\label{acoerciv}
\end{align}
cf. Hansbo and Larson \cite{HaLa11}, and, by choosing
\[
 \gamma < \frac{1}{2\vert\bfkappa^h\vert_\text{max}^2(\mu+\lambda_0)}
\]
we regain (\ref{stabdelta}). Note that if we compute $\gamma$ locally, the local $\vert\bfkappa^h\vert$ can be used instead of 
$\vert\bfkappa^h\vert_\text{max}$. Assuming we have a discrete Korn's
inequality, $\|\bfu^h\|_{H^1(\Gamma_h)}
\leq C_K a(\bfu_h,\bfu_h)$, we may introduce the norm
\[
\tnorm{\bfu^h}^2:=  a_{\Gamma_h}(\bfu^h,\bfu^h) +\gamma^{-1}\|[P_\gamma(\bfu^h)- \psi]_+\|^2_{\Gamma_h}.
\]
Observe that above we also assumed a Poincar\'e inequality. In
the numerical examples below, Korn and Poincar\'e are made to hold using midline symmetry assumptions.
Using Korn's inequality we have the boundedness of the right hand side, 
\[
|L(\bfu^h)|
\leq C_{K} \|\bff^e\|_{L^2(\Gamma_h)}\tnorm{\bfu^h}.
\]
Using \eqref{stab_rel_1}, \eqref{acoerciv}, the Cauchy-Schwarz 
inequality and the boundedness of $L$, we then obtain 
\[
C_\gamma \tnorm{\bfu^h}^2 \leq (C_{K} \|\bff^e\|_{L^2(\Gamma_h)} +
\gamma^{-\frac12} \|\psi\|_{L^2(\Gamma_h)} + \gamma C_{\kappa\mu\lambda} \|f^h_n\|_{L^2(\Gamma_h)}) \tnorm{\bfu^h},
\]
with $C_{\kappa\mu\lambda} = O(\vert\bfkappa^h\vert_\text{max}(\mu+\lambda_0)^{\frac12})$,
and consequently using the triangle inequality and the bound on $\gamma$,
\begin{equation}\label{stability_bound}
C_\gamma \tnorm{\bfu^h} \leq  (C_{K} \|\bff^e\|_{L^2(\Gamma_h)} +
\gamma^{-\frac12} \|g\|_{L^2(\Gamma_h)} +C \gamma^{\frac12} \|f^h_n\|_{L^2(\Gamma_h)}).
\end{equation}
 The existence of a unique solution can then be shown using the
 stability bound \eqref{stability_bound}, the continuity and
 the monotonicity of $b_{\Gamma_h}$ following the arguments in Burman et al. \cite{BuHaLaSt17}.

\section{Numerical examples}\label{numex}

In the numerical examples we use a $P^1$--continuous approximations of the displacements and a superparametric $P^2$--continuous approximation of the geometry.
We remark that a $P^1$--continuous geometry leads to zero curvature in each element, so a post-processing step would then be required to approximate the curvature. 
With a piecewise $P^k$ geometry approximation, for $k\geq 2$, we compute an approximate curvature directly as discussed in Section \ref{GLS}.

To visualize the reaction force, we first define
\[
p^* :=\gamma^{-1}[u^h_n- g+\gamma (f_n^h+ \sigma_\kappa(\bfu^h))]_+
\]
which is not continuous; we then perform a lumped mass $L_2$--projection of $p^*$ onto the space $W^h_1$ to obtain a smoothly varying reaction force $p^h$.

In the second numerical example, we also compare or results with the multiplier method obtained by using $P^1$--continuous approximations for both $\bfu$ and $p$, which, in our experience, leads to a stable solution.

\subsection{A sphere contained in an ellipsoid}

In this example, a sphere of radius $R=3/4$ is placed in an ellipsoid with varying fixed major axis, of length $R_\text{max}=3/2$ m, and varying minor axis, of length $R_\text{min}$. 
The parameters are as follows: Young's modulus $E=100$ MPa,
Poisson's ratio $\nu=0.5$, $\gamma =10^{-2}/(\lambda_0+\mu)$. No external load is applied.

\begin{figure}[h]
\begin{center}
\includegraphics[height=10cm]{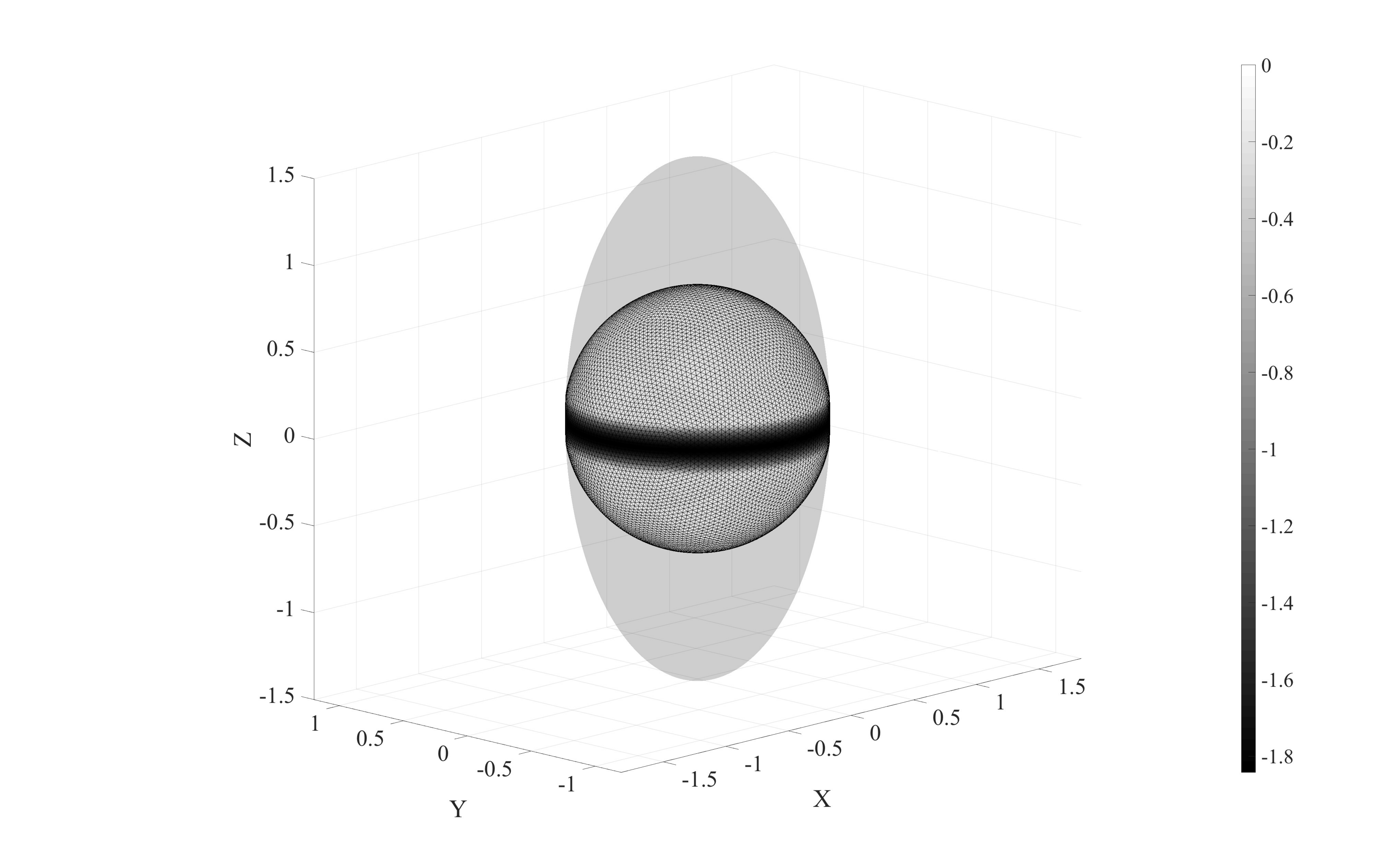}
\end{center}
\caption{Contact solution, $R_\text{min}=0.74$.\label{disp1}}
\end{figure}

\begin{figure}[h]
\begin{center}
\includegraphics[height=10cm]{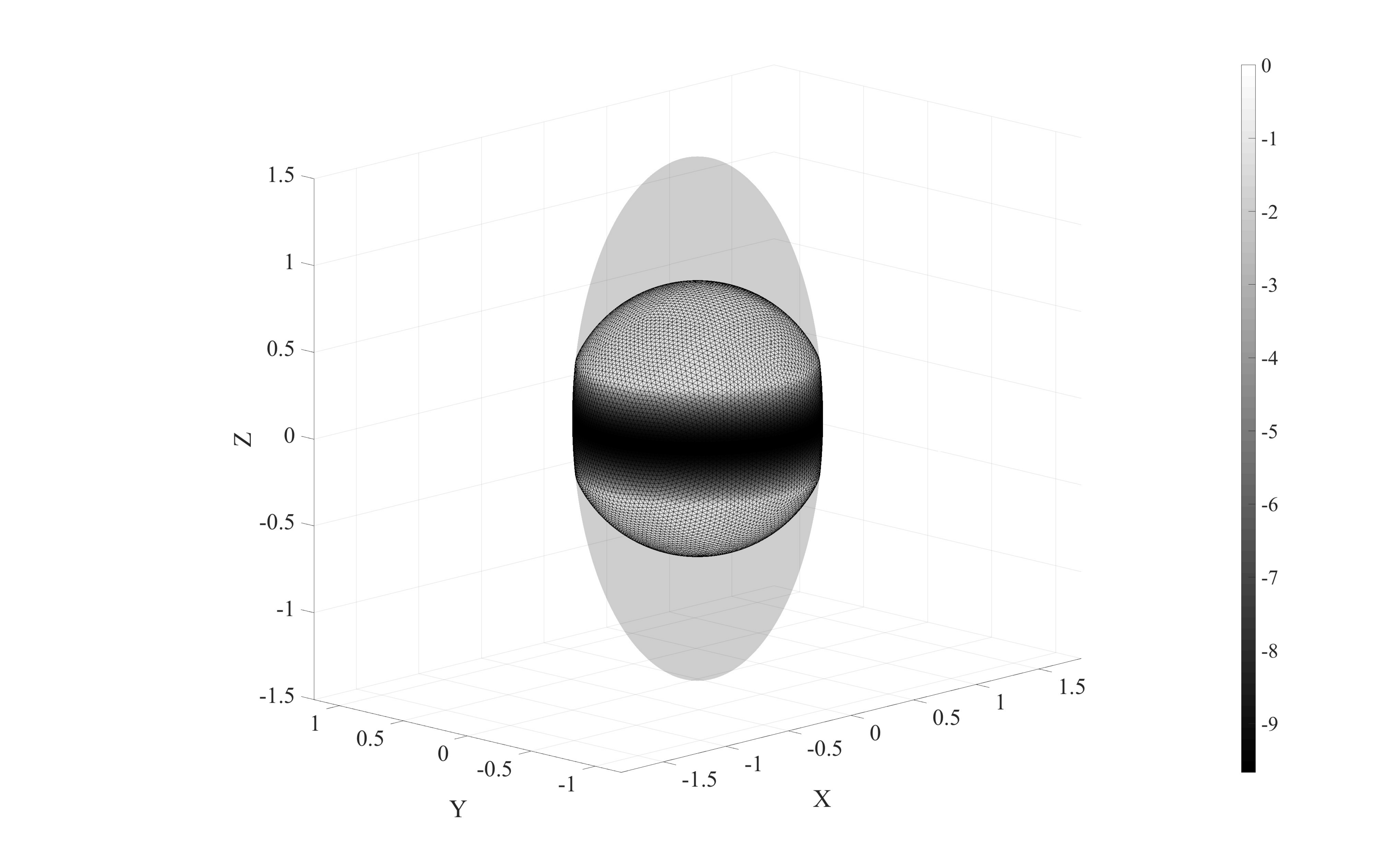}
\end{center}
\caption{Contact solution, $R_\text{min}=0.7$.\label{disp2}}
\end{figure}

\begin{figure}[h]
\begin{center}
\includegraphics[height=10cm]{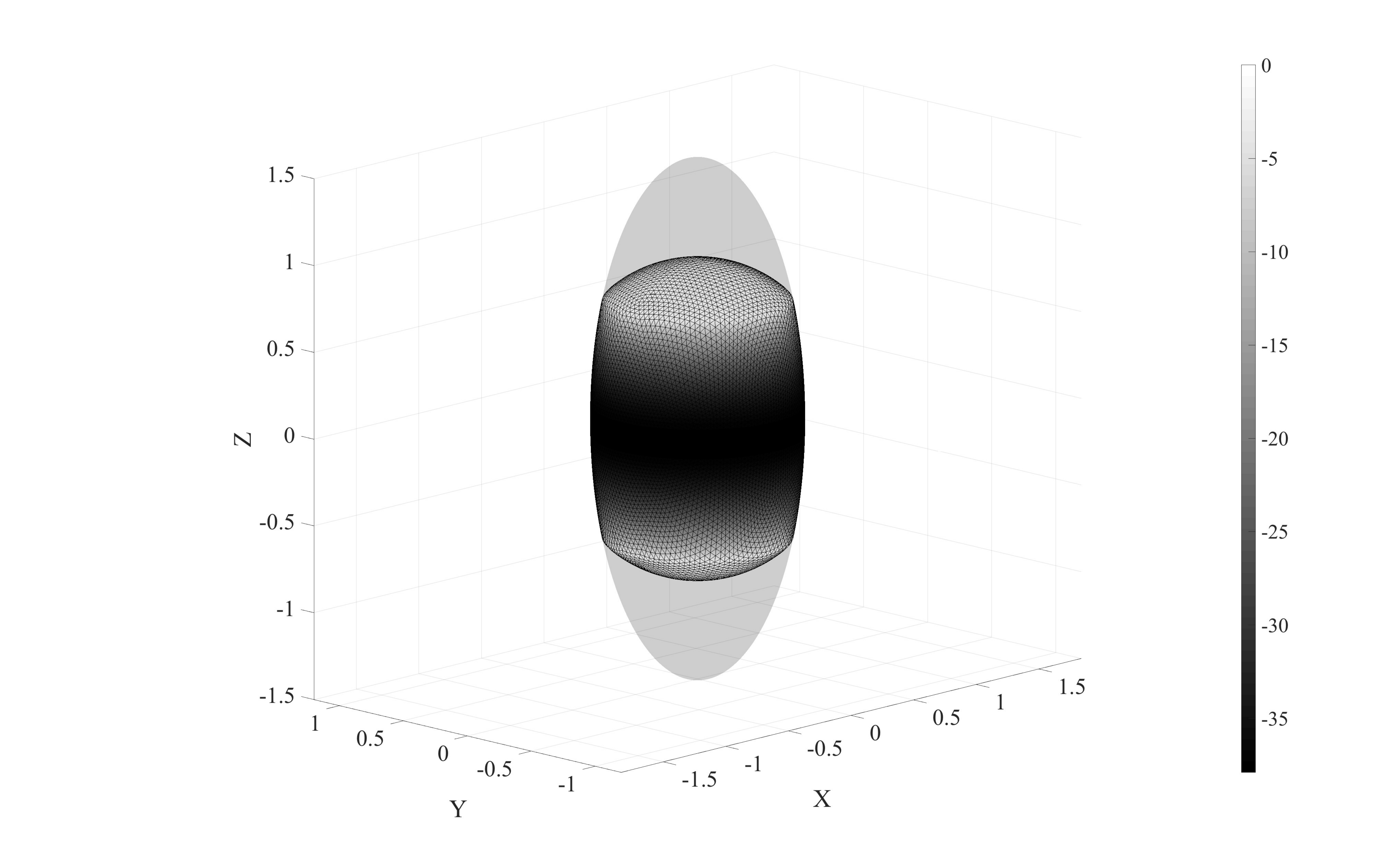}
\end{center}
\caption{Contact solution, $R_\text{min}=0.6$.\label{disp3}}
\end{figure}

In Figs. \ref{disp1}--\ref{disp3} we show the computed result on a mesh consisting of 20480 triangles and 10242 nodes with $R_\text{min}\in \{0.74,0.7,0.6\}$. We show a deformation plot of nodal displacements $(\bfu^h\cdot\bfn)\bfn$ together with isoplots of the computed reaction force. The enclosing ellipsoid is shaded.

\subsection{A sphere in contact with a rigid floor}

We use the same sphere, data, and mesh as in the previous example, now in contact with a floor located at $z=-0.74$ m and with a load $\bff= (0,0,-1)$ MPa/m$^3$.

\begin{figure}[h]
\begin{center}
\includegraphics[height=10cm]{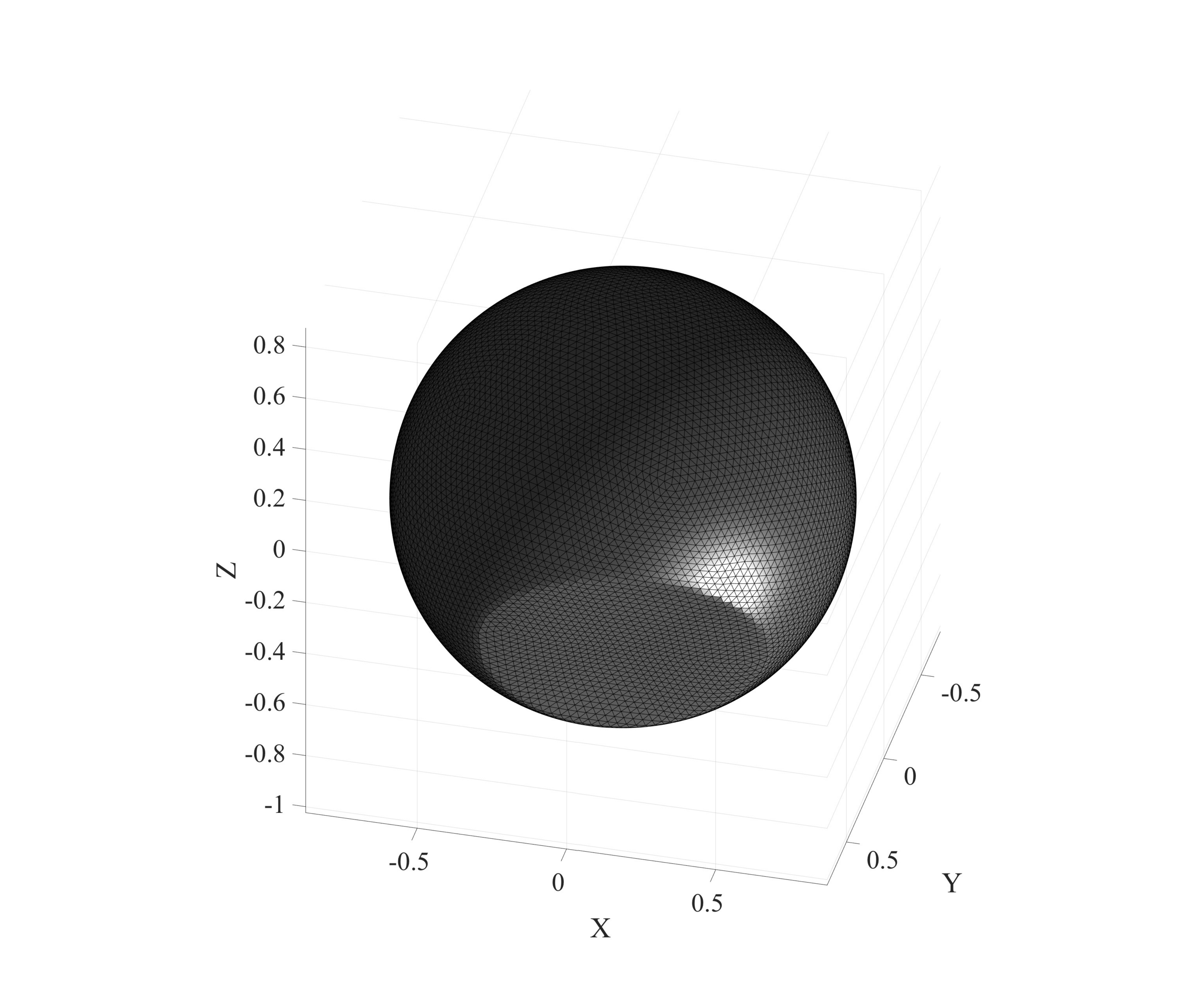}
\end{center}
\caption{Contact solution, contact with a rigid floor.\label{elevfloor}}
\end{figure}

\begin{figure}[h]
\begin{center}
\includegraphics[height=6cm]{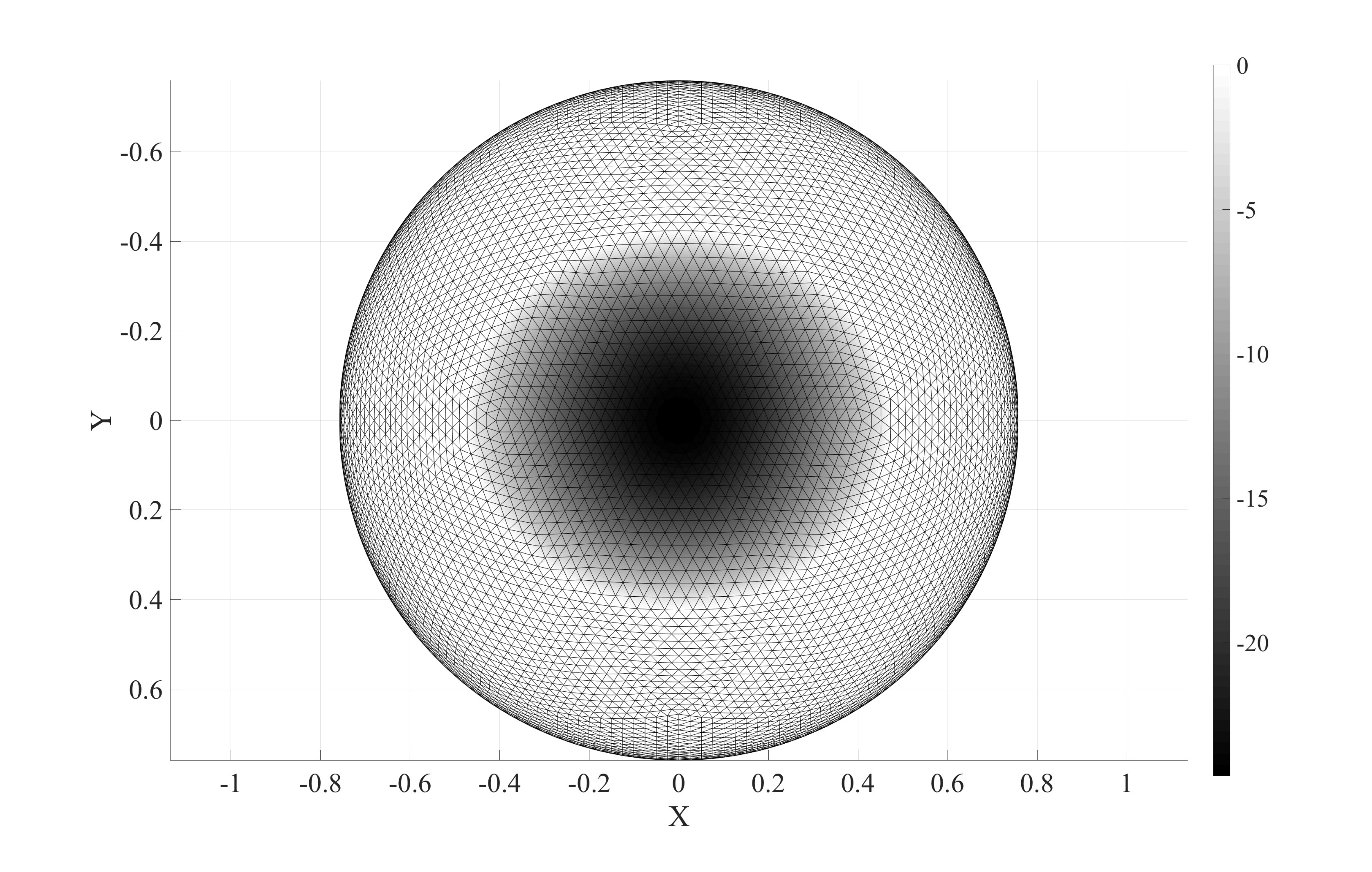}
\end{center}
\caption{Isoplot of the reaction force.\label{isofloor}}
\end{figure}

In Fig. \ref{elevfloor} we show the computed solution , and in Fig. \ref{isofloor} we give an isoplot of the reaction force.

\begin{figure}[h]
\begin{center}
\includegraphics[height=5cm]{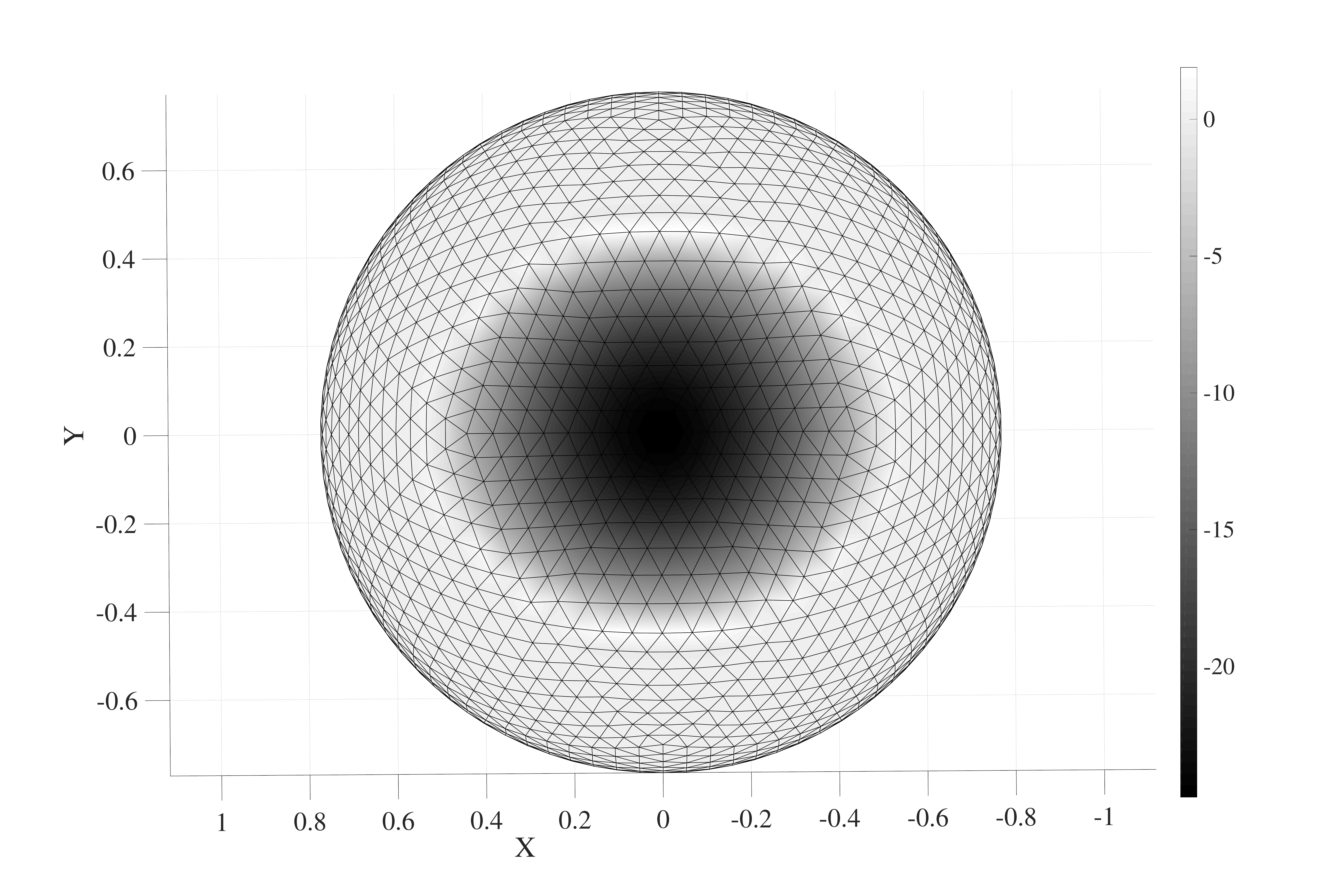}\includegraphics[height=4.9cm]{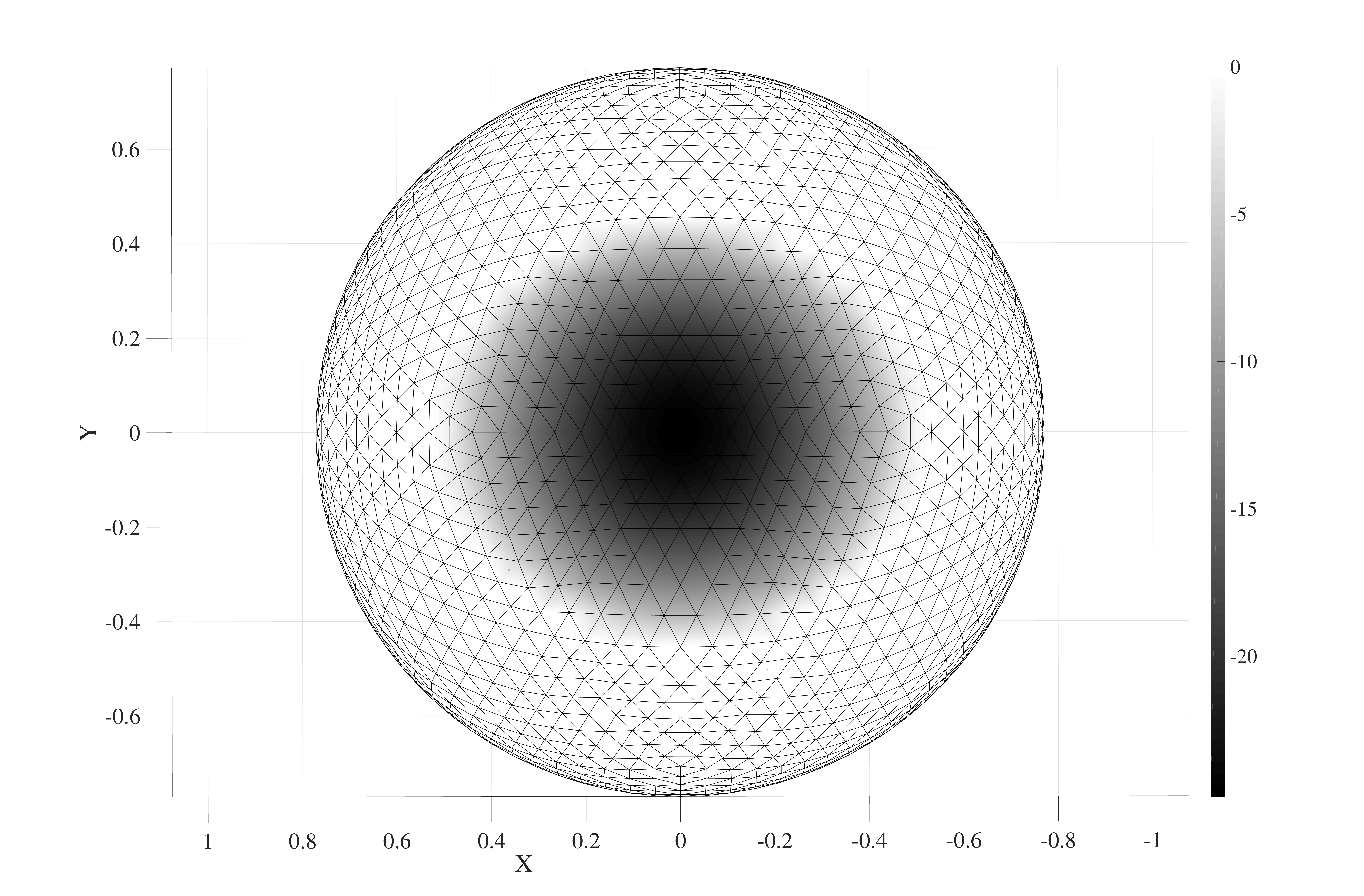}
\end{center}
\caption{Isoplot of the reaction force using a separate approximation of $p$; $p^h$ to the left and $p^*:=-(1/\gamma)(\bfn\cdot\bfu^h-g-\gamma p^h)$ to the right.\label{isofloor2}}
\end{figure}

For comparison, we also show a computation performed with $C^0$--continuous, piecewise $P^1$ contact forces. Here we the same $\gamma$ as in the GLS case and show, in Fig. \ref{isofloor2},
isoplots of the computed contact force as well as a post--processed contact force. The solution agrees with the GLS result.

\section{Concluding remarks}\label{conclusions}

We have proposed a multiplier method for the  analysis of friction free contact between curved membranes and rigid obstacles. 
By use of  a Galerkin/least squares approach we 
also show how to eliminate the multiplier, which avoids the question of {\it inf--sup}\/ stability of the combination of 
approximations for the multiplier and primary variable
and leads to a symmetric positive definite discrete system. 

\subsection*{Acknowledgements}
This research was supported in part by 
the Swedish Foundation for Strategic Research Grant No.\ AM13-0029, 
the Swedish Research Council Grant No. 2013-4708, and 
the Swedish strategic research programme eSSENCE. The first author was
supported by EPSRC Grant EP/P01576X/1.


\newpage

\begin{thebibliography}{10}

\bibitem{AlCu91}
P.~Alart and A.~Curnier.
\newblock A mixed formulation for frictional contact problems prone to {N}ewton
  like solution methods.
\newblock {\em Comput. Methods Appl. Mech. Engrg.}, 92(3):353--375, 1991.

\bibitem{BaHu91}
H.~Barbosa and T.~Hughes.
\newblock The finite element method with {L}agrange multipliers on the
  boundary: circumventing the {B}abu\v ska-{B}rezzi condition.
\newblock {\em Comput. Methods Appl. Mech. Engrg.}, 85(1):109--128, 1991.

\bibitem{BuHaLaSt17}
E.~Burman, P.~Hansbo, M.~Larson, and R.~Stenberg.
\newblock Galerkin least squares finite element method for the obstacle
  problem.
\newblock {\em Comput. Methods Appl. Mech. Engrg.}, 313:362--374, 2017.

\bibitem{ChHi12}
F.~Chouly and P.~Hild.
\newblock A {N}itsche-based method for unilateral contact problems: numerical
  analysis.
\newblock {\em SIAM J. Numer. Anal.}, 51(2):1295--1307, 2013.

\bibitem{CiLo96}
P.~Ciarlet and V.~Lods.
\newblock Asymptotic analysis of linearly elastic shells. {I}. {J}ustification
  of membrane shell equations.
\newblock {\em Arch. Rational Mech. Anal.}, 136(2):119--161, 1996.

\bibitem{CiSa93}
P.~Ciarlet and E.~Sanchez-Palencia.
\newblock Un th\'eor\`eme d'existence et d'unicit\'e pour les \'equations des
  coques membranaires.
\newblock {\em C. R. Acad. Sci. Paris S\'er. I Math.}, 317(8):801--805, 1993.

\bibitem{DeZo97}
M.~Delfour and J.-P. Zol\'esio.
\newblock Differential equations for linear shells: comparison between
  intrinsic and classical models.
\newblock In {\em Advances in mathematical sciences: {CRM}'s 25 years
  ({M}ontreal, {PQ}, 1994)}, volume~11 of {\em CRM Proc. Lecture Notes}, pages
  41--124. Amer. Math. Soc., Providence, RI, 1997.

\bibitem{FG83}
M.~Fortin and R.~Glowinski.
\newblock {\em Augmented {L}agrangian {M}ethods}.
\newblock North-Holland Publishing Co., Amsterdam, 1983.

\bibitem{GuMo75}
M.~Gurtin and A.~Murdoch.
\newblock A continuum theory of elastic material surfaces.
\newblock {\em Arch. Rational Mech. Anal.}, 57:291--323, 1975.

\bibitem{HaLa11}
P.~Hansbo and M.~Larson.
\newblock A posteriori error estimates for continuous/discontinuous {G}alerkin
  approximations of the {K}irchhoff-{L}ove plate.
\newblock {\em Comput. Methods Appl. Mech. Engrg.}, 200(47-48):3289--3295,
  2011.

\bibitem{HaLa14}
P.~Hansbo and M.~Larson.
\newblock Finite element modeling of a linear membrane shell problem using
  tangential differential calculus.
\newblock {\em Comput. Methods Appl. Mech. Engrg.}, 270:1--14, 2014.

\bibitem{HaLaLa15}
P.~Hansbo, M.~Larson, and F.~Larsson.
\newblock Tangential differential calculus and the finite element modeling of a
  large deformation elastic membrane problem.
\newblock {\em Comput. Mech.}, 56(1):87--95, 2015.

\bibitem{HeHa06}
P.~Heintz and P.~Hansbo.
\newblock Stabilized {L}agrange multiplier methods for bilateral elastic
  contact with friction.
\newblock {\em Comput. Methods Appl. Mech. Engrg.}, 195(33-36):4323--4333,
  2006.

\bibitem{HiRe10}
P.~Hild and Y.~Renard.
\newblock A stabilized {L}agrange multiplier method for the finite element
  approximation of contact problems in elastostatics.
\newblock {\em Numer. Math.}, 115(1):101--129, 2010.

\bibitem{Kr91}
E.~Kreyszig.
\newblock {\em Differential geometry}.
\newblock Dover Publications Inc., New York, 1991.

\bibitem{OlHa09}
J.~Oliver, S.~Hartmann, J.~Cante, R.~Weyler, and J.~Hern\'andez.
\newblock A contact domain method for large deformation frictional contact
  problems. {I}. {T}heoretical basis.
\newblock {\em Comput. Methods Appl. Mech. Engrg.}, 198(33-36):2591--2606,
  2009.

\end{thebibliography}
\end{document}